\documentstyle[amstex,amssymb,12pt]{amsart}

\input xypic.tex

\theoremstyle{plain}
\newtheorem{thm}{Theorem}[section]
\newtheorem{lem}[thm]{Lemma}
\newtheorem{prop}[thm]{Proposition}
\newtheorem{cor}[thm]{Corollary}

\theoremstyle{definition}
\newtheorem{alg}[thm]{Algorithm}
\newtheorem{df}[thm]{Definition}

\newtheorem{ex}[thm]{Example}

\newtheorem{rem}[thm]{Remark}

\def\ann{\operatorname{ann}}
\def\cd{\operatorname{cd}}
\def\codim{\operatorname{codim}}

\def\depth{\operatorname{depth}}
\def\ext{\operatorname{Ext}}
\def\hom{\operatorname{Hom}}
\def\height{\operatorname{ht}}
\def\im{\operatorname{im}}
\def\ini{\operatorname{in}}
\def\lcm{\operatorname{lcm}}
\def\soc{\operatorname{soc}}
\def\spec{\operatorname{Spec}}

\def\A{{\Bbb A}}
\def\C{{\Bbb C}}
\def\del{\partial}
\def\eps{\varepsilon}
\def\m{{\frak m}}
\def\N{{\Bbb N}}
\def\O{{\cal O}}
\def\P{{\Bbb P}}

\def\Z{{\Bbb Z}}

\newcommand{\eqa}{\begin{eqnarray*}}
\newcommand{\eqe}{\end{eqnarray*}}

\renewcommand{\bar}{\overline}
\def\wo{\setminus}
\def\qed{\hfill$\Box$}
\def\_#1{\underline{#1}}
\def\curlN{{\cal N}}

\def\sbe{\subseteq}
\def\into{\hookrightarrow}

\def\m{{\frak m}}
\def\qed{\hfill$\Box$}

\def\mylabel#1{\label{#1}}

\numberwithin{equation}{section}

\begin{document}

\title[Algorithms in Local Cohomology]
       {Algorithmic Computation of Local Cohomology Modules and the
       Local Cohomological Dimension of Algebraic Varieties}
\author[Uli Walther]{Uli Walther\\ University of Minnesota}
\address{School of Mathematics, University of Minnesota, Minneapolis, MN 55455}
\email{walther@@math.umn.edu}
\begin{abstract}
In this paper we present algorithms that compute certain local
cohomology modules associated to a ring of polynomials containing the
rational numbers. In particular we are able to compute the local
cohomological dimension of algebraic varieties in characteristic
zero. Our approach is based on the theory of $D$-modules.
\end{abstract}

\maketitle

\section{Introduction}
\subsection{} Let $R$ be a commutative Noetherian ring, $I$ an ideal in
$R$ and $M$ an 
$R$-module. The $i$-th {\em local
cohomology functor} with respect to $I$ is the $i$-th right
derived functor of the functor $H^0_I(-)$ which sends $M$ to the
$I$-torsion $\bigcup_{k=1}^\infty (0:_MI^k)$ of $M$ and is denoted by
$H^i_I(-)$. Local 
cohomology was introduced by Grothendieck as an algebraic analog of
(classical) relative cohomology. A brief introduction to local
cohomology may be found in appendix 4 of \cite{E}. 

The {\em
cohomological dimension} of $I$ in $R$, denoted by $\cd(R,I)$, is the
smallest integer $c$ such that the local cohomology modules
$H^q_I(M)=0$ for all $q>c$ and all $R$-modules $M$. If $R$ is the
coordinate ring of an affine variety $X$ and $I\sbe R$ is the defining
ideal of the Zariski closed subset $V\sbe X$ then the {\em local
cohomological dimension} of $V$ in $X$ is defined as $\cd(R,I)$.
It is not hard to show that if $X$ is smooth, then the integer
$\dim(X)-\cd(R,I)$ depends only on $V$ but neither on $X$ nor on the
embedding $V\into X$. 
\subsection{} Knowledge of local cohomology modules provides interesting
information, illustrated by the following three situations. Let $I\sbe
R$ and $c=\cd(R,I)$. Then 
$I$ cannot be generated by fewer than $c$ elements. In fact, no ideal $J$ 
with the same radical as $I$ will be generated by fewer than $c$ elements.

Let $H^i_{dR}$ stand for the $i$-th de Rham cohomology group. 
A second application is a family of results commonly known as Barth
theorems which are a generalization of the classical
Lefschetz theorem that states that if $Y\sbe \P^n_\C$ is a hypersurface 
then $H^i_{dR}(\P^n_\C)\to
H^i_{dR}(Y)$ is an isomorphism for $i<\dim (Y)-1$ and injective for
$i=\dim (Y)$. For example, let $Y\sbe \P_\C^n$ be a closed subset and $I\sbe
R=\C[x_0,\ldots,x_n]$ the 
defining ideal of $Y$. 
Then $H^i_{dR}(\P^n_\C)\to H^i_{dR}(Y)$ is an isomorphism for $i\le
\depth_{\O_{\P^n_\C}}(\O_Y)-\cd(R,I)$ (compare \cite{Og}, 4.7 and
\cite{DRCAV}, the theorem after III.7.6).

Finally, it is also a consequence of the work of Ogus and Hartshorne
(\cite{Og}, 2.2, 2.3 and \cite{DRCAV}, IV.3.1) that if $I\sbe
R=\C[x_0,\ldots,x_n]$ is the defining 
ideal of a complex smooth variety $V\sbe \P^n_\C$ then, for $i<n-\codim (V)$,
$\dim_\C\soc_R
(H^0_\m(H^{n-i} (R)))$ equals $\dim_\C H^i_x(\tilde V,\C)$ where $H^i_x(\tilde
V,\C)$ stands for the $i$-th singular cohomology group of the affine
cone $\tilde V$ over $V$ with support in the vertex $x$ of $\tilde V$ and
with coefficients in $\C$ ($\soc_R(M)$  denotes the socle
$(0:_M\m)\sbe M$ for any 
$R$-module $M$).
\subsection{} 
The cohomological dimension has been studied by many authors, for
example R.~Hartshorne (\cite{CDAV}), A.~Ogus (\cite{Og}),
R.~Hartshorne and R.~Speiser
(\cite{H-Sp}), C.~Peskine and L.~Szpiro (\cite{P-S}), G.~Faltings
(\cite{F}), 
C.~Huneke and G.~Lyubeznik (\cite{Hu-L}). 
Yet 
despite this extensive effort, the  
problem of finding an algorithm for the computation of cohomological
dimension remained open. For the determination of $\cd(R,I)$ it is in
fact enough to find an algorithm to 
decide whether or not the local cohomology module $H^i_I(R)=0$ for
given $i, R, I$. This is because $H^q_I(R)=0$ for all $q>c$ implies
$\cd(R,I)\le c$ (see \cite{CDAV}, section 1).  In \cite{L-Fmod}
G.~Lyubeznik gave an
algorithm for deciding whether or not 
$H^i_I(R)=0$ for all $I\sbe R=K[x_1,\ldots,x_n]$ where
$K$ is a field of positive characteristic.
One
of the main purposes of 
this work is to produce such an algorithm in the case where $K$ is a field
containing the rational numbers and $R=K[x_1,\ldots,x_n]$.

Since in such a situation
the local cohomology
modules $H^i_I(R)$ have a natural
structure of finitely generated left $D(R,K)$-modules (\cite{L-Dmod}), $D(R,K)$
being the ring 
of $K$-linear differential operators of $R$, explicit computations may be 
performed.
Using 
this finiteness we employ the theory of
Gr\"obner bases
to develop 
algorithms that give a representation of $H^i_I(R)$ and $H^i_\m (H^j_I(R))$
for all triples  $i,j\in \N,I\subseteq R$ in terms of generators and
relations over 
$D(R,K)$ (where $\m=(x_1,\ldots,x_n)$). This also leads to an
algorithm for the computation of the invariants 
$\lambda_{i,j}(R/I)=\dim_K\soc_R(H^i_\m(H^{n-j}_I(R)))$ introduced in
\cite{L-Dmod}. 

We remark that if $R$ is an arbitrary finitely generated $K$-algebra
and $I$ is an ideal in $R$
then, if 
$R$ is
regular, our algorithms can be used to determine $\cd(R,I)$ for
all ideals $I$ of $R$, 
but if $R$ is not regular, then the problem of algorithmic
determination of $\cd(R,I)$ remains open (see also the comments in
subsection \ref{singular_spaces}). 
\subsection{}
The outline of the paper is as follows. 
The next section is devoted to a short survey of results on
local cohomology and $D$-modules as they apply to our work, as well as
their interrelationship. 

In section \ref{sec-gb} we
review the theory of Gr\"obner bases as it applies to $A_n$ and
modules over the Weyl algebra. Most of that section should be
standard and readers interested in proofs and more details are
encouraged to look at the book by D.~Eisenbud (\cite{E}, chapter 15
for the commutative case)
or the fundamental article \cite{KR-W} (for the more general situation. 

In section \ref{sec-mal-kash} we generalize some results due to
B.~Malgrange and 
M.~Kashiwara on $D$-modules and their localizations. The purpose of
sections \ref{sec-mal-kash} and \ref{sec-oaku} is to find a 
representation of $R_f\otimes N$ as a cyclic $A_n$-module if $N$
is a given holonomic $D$-module (for a definition and some properties
of holonomic modules, see subsection \ref{D-modules} below).
Many of the essential ideas in section \ref{sec-oaku} come from
T.~Oaku's work \cite{Oa}. 

In
section \ref{sec-lc} we describe our main results, namely 
algorithms that for
arbitrary $i,j,k,I$ 
determine the structure of
$H^k_I(R), H^i_\m (H^j_I(R))$ and find $\lambda_{i,j}(R/I)$. 
Some of these algorithms have been 
implemented in the 
programming language C and the theory is illustrated with examples.
The final section is devoted to comments on 
implementations, effectivity and examples.

It is a pleasure to thank my advisor Gennady Lyubeznik for suggesting the
problem of algorithmic computation of cohomological dimension 
to me and pointing out that the theory of $D$-modules might be
useful for its solution.

\section{Preliminaries}
\label{sec-prelim}
\subsection{Notation} Throughout we shall use the following notation: $K$ will
denote a field 
of characteristic zero, $R=K[x_1,\ldots,x_n]$ the ring of polynomials
over $K$ in $n$ variables, $A_n=K\langle
x_1,\del_1,\ldots,x_n,\del_n\rangle$ the Weyl algebra over $K$ in $n$
variables, or, equivalently, the ring of $K$-linear differential operators on
$R$, $\m$ will stand for the maximal ideal $(x_1,\ldots,x_n)$ of
$R$, $\Delta$ will denote the maximal left ideal $(\del_1,\ldots,\del_n)$
of $A_n$ and  $I$ will stand for the ideal $(f_1,\ldots,f_r)$ in $R$. 

All tensor products in this work will be over $R$ and all
$A_n$-modules (resp.~ideals) will be left modules (resp.~left ideals).
\subsection{Local cohomology} 
It turns out that
$H^k_I(M)$ may be 
computed as follows. Let $C^\bullet(f_i)$ be the complex $0\to
R\stackrel{1\to\frac{1}{1}}{\longrightarrow} R_{f_i}\to 0$. Then
$H^k_I(M)$ is the $k$-th 
cohomology group of the {\em \v Cech complex} defined by 
$C^\bullet(M;f_1,\ldots,f_r) 
=\bigotimes_1^r C^\bullet(f_i)\otimes M$. 
Unfortunately,
explicit 
calculations are complicated by the fact that $H^k_I(M)$ is rarely
finitely generated as $R$-module. This
difficulty disappears for $H^k_I(R)$ if we 
enlarge the ring to $A_n$, in essence because $R_f$ is finitely
generated over $A_n$ 
for all $f\in R$.
\subsection{$D$-modules} \mylabel{D-modules}
A good introduction to $D$-modules is the book by Bj\"ork, \cite{B}.

Let $f\in R$. Then the $R$-module $R_f$ has a
structure as left $A_n$-module: $x_i(\frac{g}{f^k})=\frac{x_ig}{f^k},
\del_i(\frac{g}{f^k})=\frac{\del_i(g)f-k\del_i(f)g}{f^{k+1}}$. This
may be thought of as a special case of localizing an $A_n$-module: if
$M$ is an $A_n$-module and $f\in R$ then $R_f\otimes_R M$ becomes an
$A_n$-module via $\del_i(\frac{g}{f^k}\otimes
m)=\del_i(\frac{g}{f^k})\otimes m+\frac{g}{f^k}\otimes
\del_i m$. Localization of $A_n$-modules lies at the heart of our arguments. 

Of particular interest are the {\em
holonomic} modules which are those finitely generated $A_n$-modules $N$
for which $\ext^j_{A_n}(N,A_n)$ vanishes 
unless $j=n$. Holonomic modules are
always cyclic and of finite length over $A_n$. Besides that, if
$N=A_n/L$, $f\in
R$, $s$ is an indeterminate and $g$ is 
some fixed generator of $N$, then
there is a nonzero polynomial $b(s)$ in $K[s]$ and an
operator $P(s)\in A_n[s]$ such that $P(s)(f\cdot f^s\otimes g)=b(s)\cdot
f^s\otimes g$. The unique monic polynomial that divides 
all other polynomials satisfying an identity of this type is called the
{\em Bernstein polynomial} of $L$ and $f$ and denoted by $b_f^L(s)$.
Any operator $P(s)$ that satisfies $P(s)f^{s+1}\otimes
g=b_f^L(s)f^s\otimes g$ we shall call a {\em Bernstein operator} and refer
to the roots of $b_f^L(s)$ as {\em Bernstein roots} of $f$ on $A_n/L$.

Localizations of 
holonomic modules at a single element (and hence at any finite number
of elements) of $R$ are holonomic (see \cite{B}, section 5.9) and 
in particular cyclic over $A_n$, generated by $f^{-a}g$ for
sufficiently large $a\in {\Bbb N}$ (see also our proposition
\ref{kashiwara}). So the complex $C^\bullet(N;f_1,\ldots,f_r)$  
consists of holonomic $A_n$-modules whenever $N$ is holonomic. 
This facilitates the use
of Gr\"obner bases as 
computational tool for maps between holonomic modules and their
localizations. 
As a
special case we note that 
localizations of $R$ are holonomic, and hence finite,  over $A_n$ (since
$R=A_n/\Delta$ is holonomic). 
\subsection{The \v Cech complex}
In \cite{L-Dmod} it is shown that local cohomology modules are not only
$D$-modules but in fact holonomic: we know already that the modules in
the \v Cech complex are holonomic, it suffices to show that the maps are
$A_n$-linear. All maps in the \v Cech complex are direct sums of
localization maps. Suppose $R_f$ is generated by $f^s$ and
$R_{fg}$ by
$(fg)^t$. We may replace $s,t$ by their minimum $u$ and then we
see that the inclusion $R_f\to R_{fg}$ is nothing but the map
$A_n/\ann(f^u)\to A_n/\ann((fg)^u)$ sending the coset of the operator
$P$ to the coset of the operator $P\cdot g^u$. So
$C^i(N;f_1,\ldots,f_r)\to C^{i+1}(N;f_1,\ldots,f_r)$ is an
$A_n$-linear map between holonomic modules for every holonomic $N$.
One can prove that kernels and cokernels of $A_n$-linear maps between 
holonomic modules are holonomic. Holonomicity of $H^k_I(R)$
follows. 

In the same way it can be seen that $H^i_\m (H^j_I(R))$ is holonomic for
$i,j\in{\Bbb N}$ (since $H^j_I(R)$ is holonomic).  

\section{Gr\"obner bases of modules over the Weyl algebra}
\mylabel{sec-gb}
In this section we review some of the concepts and results related to
the Buchberger algorithm in modules over Weyl
algebras. It turns out that with a little care many of the important
constructions from the theory of commutative Gr\"obner bases carry
over to our case. For an introduction into non-commutative monomial
orders and related 
topics, \cite{KR-W} is a good source.

Let us agree that every time we write an element in $A_n$, we
write it as a sum of terms $c_{\alpha\beta}x^\alpha \del^\beta$ in multi-index
notation. That is, $\alpha,\beta\in \N^n$, $c_{\alpha\beta}$ are scalars,
$x^\alpha=x_1^{\alpha_1}\cdot \ldots\cdot x_n^{\alpha_n},
\del^\beta=\del_1^{\beta_1}\cdot\ldots\cdot \del_n^{\beta_n}$ and in
every monomial 
we write first the powers of $x$ 
and then the powers of the differentials. Further, if
$m=c_{\alpha\beta}x^\alpha\del^\beta, c_{\alpha\beta}\in K$, we will
say that $m$ has degree $\deg 
m=|\alpha+\beta|$ and an operator $P\in A_n$ has degree equal to the largest
degree of any monomial occuring in $P$.

Recall that a {\em monomial order} $<$ in $A_n$ is a total order on the
monomials of $A_n$, subject to $m<m'\Rightarrow mm''<m'm''$ for all nonzero
monomials $m,m',m''$. 
Since the product of two monomials in our notation is not likely to be
a monomial (as $\del_i x_i=x_i\del_i+1$) it is not obvious that such
orderings exist at all. However, the commutator of any two monomials
$m_1,m_2$ 
will be a polynomial of degree at most $\deg m_1+\deg m_2-2$. That means
that the degree of an operator and its component of maximal degree 
is independent of the way it is
represented. Thus we may for example introduce a monomial order on
$A_n$ by taking 
any monomial order on $\tilde A_n=K[x_1,\ldots,x_n,\del_1,\ldots,\del_n]$ (the
polynomial ring in $2n$ variables) that refines the partial order
given by total degree, and saying that $m_1>m_2$ in $A_n$
if and only if $m_1>m_2$ in $\tilde A_n$. 

Let $<$ be a monomial order on $A_n$.
Let $G=\bigoplus_1^dA_n\cdot \gamma_i$ be the free
$A_n$-module on the symbols $\gamma_1,\ldots,\gamma_d$. We define a
monomial order on $G$ by $m_i\gamma_i>m_j\gamma_j$ if either $m_i>m_j$
in the order on $A_n$, or $m_i=m_j$ and $i>j$.
The largest monomial $m\gamma$ in an element $g\in G$ will be denoted by
$\ini(g)$. Of fundamental importance is 
\begin{alg}[Remainder]
\mylabel{remainder}
Let $h$ and $\_ g=\{g_i\}_1^s$ be elements of $G$. Set $h_0=h, \sigma_0=0,
j=0$ and let $\eps_i=\eps(g_i)$ be symbols. Then 
\[
\begin{array}{llll}
&{\tt Repeat}&&\\
&&{\tt If} \ini(g_i)|\ini(h_j) {\tt \,set}&\\ 
&&&\{h_{j+1}:=h_j-\frac{\ini(h_j)}{\ini(g_i)}g_i,\\
&&&\sigma_{j+1}:=\sigma_j+\frac{\ini(h_j)}{\ini(g_i)}\eps_i,\\
&&&j:=j+1\}\\
&{\tt Until}&{\tt No} \ini(g_i)|\ini(h_j).
\end{array}
\]
The result is $h_a$, called a {\em remainder $\Re(h,\_ g)$ of $h$
under division by $\_ g$}, and an expression $\sigma_a=\sum_{i=1}^s
a_i\eps_i$ 
with $a_i\in A_n$. 
By Dickson's lemma (\cite{KR-W}, 1.1) the
algorithm terminates. It is worth mentioning that $\Re(h,\_ g)$ is not
uniquely determined, it depends on which $g_i$ we pick amongst those
whose initial term divides the initial term of $h_j$.

Note that if $h_a$ is zero, $\sigma_a$ tells us how to write $h$ in
terms of $\_ g$. Such a $\sigma_a$ is called a \em{standard
expression for $h$} with respect to $\{g_1,\ldots,g_s\}$. 
\end{alg}
\begin{df}\mylabel{schreyer}
If $\ini(g_i)$ and $\ini(g_j)$ involve the same basis element of $G$,
then we set $s_{ij}=m_{ji}g_i-m_{ij}g_j$ and
$\sigma_{ij}=m_{ji}\eps_i-m_{ij}\eps_j$ where
$m_{ij}=\frac{\lcm(\ini(g_j),\ini(g_i))}{\ini(g_j)}$. Otherwise,
$\sigma_{ij}$ and $s_{ij}$ are defined to be zero. $s_{ij}$ is the
{\em Schreyer-polynomial} to $g_i$ and $g_j$. 

Suppose $\Re(s_{ij},\_ g)$ is zero for all $i,j$. Then we call $\_ g$
a {\em Gr\"obner basis} for the module $A_n\cdot(g_1,\ldots,g_s)$.
\end{df}
The following proposition (\cite{KR-W}, Lemma 3.8) indicates the
usefulness of Gr\"obner 
bases.
\begin{prop}
\mylabel{gb-char}
Let $\_ g$ be a finite set of elements of $G$. Then $\_ g$ is a Gr\"obner
basis if and only if $h\in A_n\_ g$ implies $\exists i:
\ini(g_i)|\ini(h)$.\qed 
\end{prop}
Computation of Gr\"obner bases over the Weyl algebra works just as
over polynomial rings:
\begin{alg}[Buchberger]
\mylabel{buchberger}

Input: $\_ g=\{g_1,\ldots,g_s\}\sbe G$.

Output: a Gr\"obner basis for $A_n\cdot(g_1,\ldots,g_s)$.

Begin.
\[
\begin{array}{llll}
&{\tt Repeat}&&\\
&&{\tt If\,} h=\Re(s_{ij},\underline g)\not =0&\\ 
&&&{\tt add\,} h {\tt \,to \,} \_ g\\ 
&{\tt Until}&{\tt \, all\,} \Re(s_{ij},\_ g)=0.\\
&{\tt Return } \,\,\_ g.
\end{array}
\]
\indent End.
\end{alg}
\subsection{}
\mylabel{free-kernel}
Now we shall describe the construction of kernels of $A_n$-linear maps
using Gr\"obner bases. Again, this is similar to the commutative
case and we first consider the case of a map between free $A_n$-modules.

Let $E=\bigoplus_1^sA_n\eps_i, G=\bigoplus_1^rA_n\gamma_j$ and $\phi:E\to G$
be a
$A_n$-linear map. Assume $\phi(\eps_i)=g_i$. Suppose that in order to
make $\_ g$ a Gr\"obner basis we have to add $g'_1,\ldots,g'_{s'}$ to
$\underline 
g$ which 
satisfy $g'_i=\sum_{k=1}^s a_{ik}g_k$. We get an induced map 

$\diagram
\bigoplus_1^{s+s'}A_n\eps_i\dto_\pi\drto_{\tilde\phi}\\
\bigoplus_1^sA_n\eps_i\rto_\phi &\bigoplus_1^rA_n\gamma_j
\enddiagram
$ where $\pi$ is the identity on $\eps_i$ for $i\le s$ and sends
$\eps_{i+s}$ into $\sum_{k=1}^s a_{ik}\eps_k$. Of course,
$\tilde\phi=\phi\pi$.  

The kernel of $\phi$ is just the image of the kernel of $\tilde\phi$
under $\pi$. So in order to find kernels of maps between free modules
one may assume that the generators of the source are mapped to a
Gr\"obner basis and replace $\phi$ by $\tilde\phi$. Recall from
definition \ref{schreyer} that $\sigma_{ij}=m_{ji}\eps_i-m_{ij}\eps_j$ or
zero, depending on the leading terms of $g_i$ and $g_j$.

\begin{prop}
\mylabel{syz}
Assume that $\{g_1,\ldots,g_s\}$ is a Gr\"obner
basis. Let $s_{ij}=\sum d_{ijk}g_k$ be standard expressions for the Schreyer
polynomials. Then $\{\sigma_{ij}-\sum_k d_{ijk}\eps_k\}_{1\le i<j\le s}$
generate the kernel 
of $\phi:\bigoplus_1^sA_n\eps_i\to \bigoplus_1^rA_n\gamma_j$, sending
$\eps_i$ to $g_i$.
\end{prop}
The proof proceeds exactly as in the commutative case, see for example
\cite{E}, section 15.10.8.
\subsection{}
\mylabel{kernel}
We explain now how to find a set of generators for the kernel of an
arbitrary 
$A_n$-linear map. Let $E, G$ be as in subsection \ref{free-kernel} and suppose
$A_n(p_1,\ldots,p_a)=P\sbe E, 
A_n(q_1,\ldots,q_b)=Q\sbe G$ and  $\phi:\bigoplus_1^s
A_n\eps_i/P\to \bigoplus_1^r 
A_n\gamma_j/Q$. It will 
be sufficient to consider the case $P=0$ since we may lift $\phi$ to the free
module $E$ surjecting onto $E/P$. 

Let as before $\phi(\eps_i)=g_i$. A kernel element in $E$ is a sum 
$\sum_ia_i\eps_i, a_i\in A_n$, which if $\eps_i$ is replaced by
$g_i$ can be written in terms 
of the generators $q_j$ of $Q$. Let $\_
\beta=\{\beta_1,\ldots,\beta_c\}$ be such that $\_ g\cup \_ q\cup
\_\beta$ is a Gr\"obner basis for $A_n(\_ g,\_ q)$. We may assume
that the $\beta_i$ are the results of applying algorithm
\ref{buchberger} to $\_ 
g\cup\_ q$. Then from algorithm $\ref{remainder}$ we have expressions 
\begin{equation}
\label{star}
\beta_i=\sum_j c_{ij}g_j+\sum_k c'_{ik}q_k,
\end{equation}
with $c_{ij}, c'_{ik}\in A_n$. 
Furthermore, by proposition \ref{syz}, algorithm \ref{buchberger} returns a
generating set $\_\sigma$ 
for the
syzygies on $\_g\cup\_ q\cup\_\beta$. Write
\begin{equation}
\sigma_i=\sum_j a_{ij}\eps_{g_j}+\sum_k a'_{ik}\eps_{q_k}+\sum_l
a''_{il}\eps_{\beta_l}
\end{equation}
and eliminate the last sum using the relations (\ref{star}) to obtain
syzygies
\begin{equation}
\tilde \sigma_i=\sum_j a_{ij}\eps_{g_j}+\sum_k a'_{ik}\eps_{q_k}+\sum_l
a''_{il}\left(\sum_v c_{lv}\eps_{g_v}+\sum_w c'_{lw}\eps_{q_w}\right).
\end{equation}
These will then form a generating set for the syzygies on
$\_g\cup \_q$. Cutting off the $q$-part of these syzygies
we get a set of forms 
\[
\left\{\sum_j a_{ij}\eps_{g_j}+\sum_l
a''_{il}\left(\sum_v c_{lv}\eps_{g_v}\right)\right\}
\] 
which generate the kernel of the map
$E\to G/Q$.
\subsection{}
The comments in this subsection will find their application in algorithm
\ref{lclc-alg} which computes the structure of $H^i_\m(H^j_I(R))$ as
$A_n$-module. 
Let
\mylabel{double-kernel}
\[
\diagram
M_3'\rto^\alpha& M_3\rto^{\alpha '}& M_3''\\
M_2'\rto^\beta\uto_{\phi'}& M_2\rto^{\beta'}\uto_{\psi'}& M_2''\uto_{\rho'}\\
M_1'\rto^\gamma\uto_{\phi}& M_1\rto^{\gamma'}\uto_\psi& M_1''\uto_{\rho}
\enddiagram
\]
be a commutative diagram of $A_n$-modules. 
Note that the row cohomology of the column co\-ho\-mo\-lo\-gy at $N$ is
given by 
\[
\left[\ker(\psi ')\cap {\beta
'}^{-1}(\im\rho)+\im(\psi)\right]\,\,
/\,\,\left[\beta(\ker(\phi '))+\im(\psi )\right].
\]
 In
order to compute this we need to be able to find:
\begin{itemize}
\item preimages of submodules,
\item kernels of maps,
\item intersections of submodules.
\end{itemize}
It is apparent that the second and third calculation is a special case
of the first: kernels are preimages of zero, intersections are
images of preimages (if ${A_n}^r\stackrel{\phi}{\rightarrow} {A_n}^s/M
\stackrel{\psi}{\leftarrow} {A_n}^t$ is given, then $\im(\phi)\cap
\im(\psi)=\psi(\psi^{-1}(\im(\phi)))$ ).

So suppose in the situation $\phi:{A_n}^r/M\to {A_n}^s/N$,
$\psi:{A_n}^t/P\to {A_n}^s/N$ we 
want to find the preimage under $\psi$ of the image of $\phi$. We may
reduce to the case where $M$ and $P$ are zero and then lift
$\phi,\psi$ to maps into ${A_n}^s$. The elements $x$ in
$\psi^{-1}(\im 
\phi)\sbe {A_n}^t$ are exactly the elements in
$\ker({A_n}^t\stackrel{\psi}{\to} {A_n}^s/N\to {A_n}^s/(N+\im\phi))$ and this
kernel can be found according to the comments in \ref{kernel}.

\section{$D$-modules after Kashiwara and Malgrange}
\mylabel{sec-mal-kash}
The purpose of this and the following section is as follows. Given $f\in
R$ and an ideal $L\sbe A_n$ such that $A_n/L$ is holonomic  and $L$ is
$f$-saturated (i.e.~$f\cdot P\in L$ only if $P\in L$), we want
to  compute the structure of the module 
$R_f\otimes A_n/L$. It turns out that it is useful to know the ideal
$J^L(f^s)$ which consists of the operators $P(s)\in A_n[s]$ that
annihilate $f^s\otimes \bar 1\in M:=R_f[s]f^s\otimes A_n/L$  where the bar
denotes cosets in $A_n/L$.
In order to find $J^L(f^s)$, we will consider the module 
$M$ over the ring $A_{n+1}=A_n\langle
t,\del_t\rangle$. It will turn out in \ref{malgrange} that one can
easily compute the ideal $J^L_{n+1}(f^s)\sbe A_{n+1}$ consisting of all
operators that kill $f^s\otimes\bar 1$. In section \ref{sec-oaku}
we will then show how to compute $J^L(f^s)$ from $J_{n+1}^L(f^s)$. 

The second basic fact in this section (proposition \ref{kashiwara})
shows how to 
compute the structure of $R_f\otimes A_n/L$ as $A_n$-module once
$J^L(f^s)$ is known.
\subsection{}
Consider $A_{n+1}=A_n\langle t,\del_t\rangle$, the Weyl algebra in
$x_1,\ldots,x_n$ and the new variable $t$. B.~Malgrange has defined an
action of $t$ 
and $\del_t$ on $M=R_f[s]\cdot f^s\otimes_R A_n/L$ by $t(g(x,s)\cdot
f^s\otimes \bar P)=g(x,s+1)f\cdot f^s\otimes \bar P$ and $\del_t(g(x,s)\cdot
f^s\otimes \bar P)=\frac{-s}{f}g(x,s-1)\cdot f^s\otimes \bar P$ for
$\bar P\in A_n/L$. $A_n$ acts on $M$ as expected, the variables by
multiplication 
on the left, the $\del_i$ by the product rule. 

One checks that this actually
defines an structure of $M$ as a left  $A_{n+1}$-module and that
$-\del_tt$ acts as 
multiplication by $s$. 

We denote by $J^L_{n+1}(f^s)$ the ideal in $A_{n+1}$ that annihilates the
element $f^s\otimes \bar 1$ in $M$. Then we have an induced morphism
of $A_{n+1}$-modules $A/J^L_{n+1}(f^s)\to M$ sending
$P+J^L_{n+1}(f^s)$ to 
$P(f^s\otimes\bar 1)$. 

The operators $t$ and $\del_t$ were introduced in \cite{M}. The
following lemma generalizes lemma 4.1 in \cite{M} (as well as
part of the proof given there) where the special case
$L=(\del_1,\ldots,\del_n), A_n/L=R$ is considered.

Note that $J^L_{n+1}(f^s)$ makes perfect sense even if
$L$ is not holonomic.

\begin{lem}
\mylabel{malgrange}
Suppose that $L=A_n\cdot (P_1,\ldots,P_r)$ is $f$-saturated (i.e., if
$f\cdot P\in L$, then $P\in L$).
With the above definitions, $J^L_{n+1}(f^s)$ is the
ideal generated by $f-t$ together with the images of the $P_j$ under
the automorphism $\phi$ of $A_{n+1}$ induced by $x\to x, t\to t-f$. 
\pf
The automorphism sends $\del_i$ to $\del_i+f_i\del_t$ and $\del_t$ to
$\del_t$. So if we write $P_j=P_j(\del_1,\ldots,\del_n)$, then $\phi
P_j=P_j(\del_1+f_1\del_t,\ldots,\del_n+f_n\del_t)$.

One checks that $(\del_i+f_i\del_t)(f^s\otimes \bar Q)=f^s\otimes \bar{\del_i
Q}$ for all differential operators $Q$ so that
$\phi(P_j(\del_1,\ldots,\del_n))(f^s\otimes 
\bar 1)=f^s\otimes \bar{P_j(\del_1,\ldots,\del_n)}=0$. By
definition, $f\cdot 
(f^s\otimes \bar 1)=t\cdot (f^s\otimes \bar 1)$. So $t-f\in
J^L_{n+1}(f^s)$ and $\phi(P_j)\in J^L_{n+1}(f^s)$ for $i=1,\ldots,r$.

Conversely let $P(f^s\otimes \bar 1)=0$. We may assume, that $P$ does not contain
any $t$ since we can eliminate $t$ using $f-t$. Now rewrite $P$ in
terms of $\del_t$ and the $\del_i+f_i\del_t$. Say, $P=\sum
c_{\alpha\beta}\del_t^\alpha
x^\beta Q_{\alpha\beta}(\del_1+f_1\del_t,\ldots,\del_n+f_n\del_t)$,
where the $Q_{\alpha\beta}$ are polynomials in $n$ variables and 
$c_{\alpha\beta}\in K$. Application to
$f^s\otimes \bar 1$ results in $\sum \del_t^\alpha(f^s\otimes
c_{\alpha\beta}x^\beta \bar{Q_{\alpha\beta}(\del_1,\ldots,\del_n)})$. 

Let $\bar\alpha$ be the largest $\alpha\in\N$ for which there is a
nonzero $c_{\alpha\beta}$ occuring in $P=\sum
c_{\alpha\beta}\del_t^\alpha
x^\beta Q_{\alpha\beta}(\del_1+f_1\del_t,\ldots,\del_n+f_n\del_t)$.
We show that the sum of terms that contain
$\del_t^{\bar\alpha}$ is in $A_{n+1}\cdot \phi(L)$ as
follows. In 
order for $P(f^s\otimes\bar 1)$ to vanish, the sum of terms with the
highest $s$-power, namely $s^{\bar\alpha}$, must vanish, 
and so $\sum_\beta c_{\bar\alpha\beta}(-1/f)^{\bar\alpha}f^s\otimes
x^\beta Q_{\bar\alpha\beta}(\del_1,\ldots,\del_n)\in
R_ff^s\otimes L$ as $R_f$ is $R$-flat.
It follows, that $\sum_\beta c_{\bar\alpha\beta}x^\beta
Q_{\bar\alpha\beta}(\del_1,\ldots,\del_n)\in L$ ($L$ is
$f$-saturated!) and hence $\sum_\beta
\del_t^{\bar\alpha}c_{\bar\alpha\beta}
x^\beta
Q_{\bar\alpha\beta}(\del_1+f_1\del_t,\ldots,\del_n+f_n\del_t)\in
A_{n+1}\cdot \phi(L)$.

So by the first part,
$P-\sum_\beta c_{\bar\alpha\beta}\del_t^{\bar\alpha} x^\beta
Q_{\bar\alpha\beta}(\del_1+f_1\del_t,\ldots,\del_n+f_n\del_t)$ kills
$f^s\otimes \bar 1$, but is of 
smaller degree in $\del_t$ than $P$ was.

The claim follows.\qed
\end{lem}

\subsection{}
\mylabel{jt-js}
Let $J^L(f^s)$ stand for the ideal in $A_n[s]\cong A_n[-\del_tt]$
that kills $f^s\otimes\bar 1\in R_f[s]f^s\otimes_R A_n/L$. Note that
$J^L(f^s)=J^L_{n+1}(f^s)\cap A_n[-\del_tt]$. Again, we may talk about
$J^L(f^s)$ independently of the holonomicity of $L$.

We will in the next section show how the lemma can be used to
determine $J^L(f^s)$. Now we show why $J^L(f^s)$ is useful, generalizing
\cite{K}, proposition 6.2.

Recall that the Bernstein polynomial $b^L_f(s)$ is defined to be the 
monic generator of the ideal of polynomials $b(s) \in K[s]$ for which
there exists an operator $P(s)\in A_n[s]$ such that $P(s)(f^{s+1}\otimes \bar
1)=b(s)f^s\otimes \bar 1$ (\cite{B}, chapter 1), and that $b_f^L(s)$
will exist for example if $L$ is holonomic. 

\begin{prop}
\mylabel{kashiwara}
If $L$ is holonomic and $a\in \Bbb Z$ is such that no integer root of
$b_f^L(s)$ is smaller than 
$a$, then we have isomorphisms 
\begin{equation}
R_f\otimes A_n/L\cong A_n[s]/J^L(f^s)|_{s=a}\cong A_n\cdot f^a\otimes \bar
1.
\end{equation}
\pf 
We mimick the proof given by Kashiwara, who proved the proposition for
the case $L=(\del_1,\ldots,\del_n), A_n/L=R$ (\cite{K}, proposition 6.2).

Let us first prove the second equality. Certainly $J^L(f^s)|_{s=a}$
kills $f^a\otimes\bar 1$. So we have to show that if $P(f^a\otimes\bar
1)=0$ then 
$P\in J^L(f^s)+A_n[s]\cdot(s-a)$. To that end note that
$st$ acts as $t(s-1)$ which means that $t\cdot
(A_n[s]/J^L(f^s))$ is a left $A_n[s]$-module. Identify $A_n[s]/J^L(f^s)$
with $\curlN^L_f:=A_n[s]\cdot (f^s\otimes\bar 1)$.
By
definition, $b^L_f(s)$ is the minimal polynomial for which there is
$P(s)$ with $b^L_f(s)(f^s\otimes \bar 1)=P(s)f^{s+1}=t\cdot
P(s-1)(f^s\otimes\bar 1)$. So $b^L_f(s)$ multiplies
$A_n[s]\cdot(f^s\otimes\bar 1)$ into $t\cdot A_n[s](f^s\otimes\bar
1)$ and whenever the polynomial $b(s)\in K[s]$ is relatively prime to
$b_f^L(s)$ 
its action on $\curlN_f^L/t\cdot \curlN_f^L$ is injective. 

Since by hypothesis $s-a+j$ is not a divisor of
$b^L_f(s)$ for $0<j\in \Bbb N$, 
\begin{equation}
\label{kash-eqn}
(s-a+j)\curlN^L_f\cap t\cdot \curlN^L_f\sbe (s-a+j)t\cdot
\curlN^L_f.
\end{equation}
So $(s-a+m)\curlN^L_f\cap
t^m\curlN^L_f\subseteq (s-a+m)t\curlN^L_f\cap t^m\curlN^L_f=t[(s-a+m-1)
\curlN^L_f\cap 
t^{m-1}\curlN^L_f]$ whenever $m\geq 1$.

We show now by induction on $m$ that $(s-a+m)\curlN^L_f\cap
t^m\curlN^L_f\sbe (s-a+m)t^m\curlN^L_f$ for $m\geq 1$. The claim
is clear for $m=1$ from equation (\ref{kash-eqn}). So let $m>1$. The
inductive hypothesis states that $(s-a+m-1)\curlN^L_f\cap
t^{m-1}\curlN^L_f\sbe (s-a+m-1)t^{m-1}\curlN^L_f$. The previous
paragraph shows that $(s-a+m)\curlN^L_f\cap t^m\curlN^L_f\sbe
t\left[(s-a+m-1)\curlN^L_f\cap t^{m-1}\curlN^L_f\right]$. Combining these two
facts we get 
\eqa
(s-a+m)\curlN^L_f\cap t^m\curlN^L_f&\sbe & t(s-a+m-1)t^{m-1}\curlN^L_f\\
&=&(s-a+m)t^m\curlN^L_f.
\eqe

Now if $P(s)\in A_n[s]$ is of degree $m$ in the $\del_i$ and
$P(a)(f^a\otimes \bar 1)=0$,
then
$P(s+m)\cdot f^m+J^L(f^s)\in (s-a+m)\cdot \curlN^L_f$ because we can
interprete 
$P(s+m)(f^{s+m}\otimes \bar 1)$ as a polynomial in $s+m$ with root
$a$. But then 
$P(s+m) (f^{s+m}\otimes\bar 1)=P(s+m)(f^{m}f^s\otimes\bar 1)$ is in 
\[
(s-a+m)\curlN^L_f\cap t^m\curlN^L_f\sbe
(s-a+m)t^m\curlN^L_f,
\]
implying $P(s+m)(f^{s+m}\otimes \bar 1)=(s-a+m)Q(s)(f^{s+m}\otimes\bar 1)$
for some $Q(s)\in A_n[s]$ (note 
that $J^L(f^s)$ kills $f^s\otimes\bar 1$). In other words, $P(s)-(s-a)Q(s-m)\in
J^L(f^s)$. 

For the first isomorphism we have to show that
$A_n\cdot(f^a\otimes\bar 1)=R_f\otimes A_n/L$. It suffices
to show that every term of the form $f^mf^a\otimes\bar Q$ is in the
module generated by $(f^a\otimes\bar 1)$ for all $m\in \Bbb
Z$. Furthermore, we may assume that $Q$ is a monomial in 
$\del_1,\ldots,\del_n$. 

Existence and definition of $b^L_f(s)$ provides an operator $P(s)$ with
$[b^L_f(s-1)]^{-1}P(s-1)(f^s\otimes \bar 1)=f^{-1}f^s\otimes \bar
1$. As $b_f^L(a-m)\not=0$ for 
all $0<m\in \Bbb N$ we have $f^{m}f^a\otimes\bar 1\in A_n\cdot
(f^a\otimes \bar 
1)$  for all $m$. Now let $Q$ be a monomial in $\del_1,\ldots,\del_n$
of $\del$-degree $j>0$ and assume that  
$f^mf^a\otimes \bar{Q'}\in A_n\cdot (f^a\otimes \bar 1)$ for all
$m$ and all operators $Q'$ of $\del$-degree lower than $j$. 
Then $Q=\del_i Q'$ for some $1\le i\le
n$. Fix $m\in\Z$. By assumption on $j$, for some $P'$ we have
$P'(f^a\otimes \bar 
1)=f^mf^a\otimes\bar{Q'}$. So
\begin{equation}
f^mf^a\otimes\bar Q=
\del_iP'(f^a\otimes\bar 1)-f_i\cdot(a+m)f^{m-1}f^a\otimes \bar{Q'}
\in A_n\cdot (f^a\otimes \bar 1).
\end{equation}
The claim follows by induction. 
This completes the proof of the proposition.\qed
\end{prop}
We remark that if any $a\in\Z $ satisfies the conditions of the
proposition, then so does every integer smaller than $a$.

\section{An algorithm of Oaku}
\mylabel{sec-oaku}
The purpose of this section is to review and generalize an algorithm
due to Oaku. 
In \cite{Oa} (algorithm 5.4.), Oaku
showed how to construct a generating set for $J^L(f^s)$ in the case where
$L=(\del_1,\ldots,\del_n)$. According to \ref{jt-js}, $J^L(f^s)$ is the
intersection of $J^L_{n+1}(f^s)$ with $A_n[-\del_tt]$. 
We shall explain how one may calculate $J\cap
A_n[-\del_tt]$ whenever $J\sbe A_{n+1}$ is any given ideal and as a
corollary develop an algorithm that for $f$-saturated $A_n/L$
computes $J^L(f^s)$. The proof
follows closely Oaku's argument.

On $A_{n+1}[y_1,y_2]$ define weights $w(t)=w(y_1)=1,
w(\del_t)=w(y_2)=-1, w(x_i)=w(\del_i)=0$. If $P=\sum_i P_i\in
A_{n+1}[y_1,y_2]$ and all $P_i$ are monomials, then we will write
$(P)^h$ for the operator $\sum_i P_i\cdot y_1^{d_i}$ where
$d_i=\max_j(w(P_j))-w(P_i)$ and call it the {\em $y_1$-homogenization}
of $P$.

Note that the
Buchberger algorithm preserves homogeneity 
in the following sense: if a set of generators for an ideal is given
and these generators are homogeneous with respect to the weights above,
then any new generator for the ideal constructed with the classical
Buchberger algorithm will also be homogeneous. (This is a consequence
of the facts that the $y_i$ commute with all other variables and that
$\del_t t=t\del_t+1$ is homogeneous of weight zero.) 

\begin{prop}
\mylabel{oaku}
Let $J=A_{n+1}\cdot(Q_1,\ldots,Q_r)$ and let
$y_1,y_2$ be two new 
variables. 
Let $I$ be the left ideal in $A_{n+1}[y_1]$ generated by the
$y_1$-homogenizations $(Q_i)^h$ of the $Q_i$, relative to the weight
$w$ above, and let $\tilde
I=A_{n+1}[y_1,y_2]\cdot (I,1-y_1y_2)$. Let $G$ be a Gr\"obner basis
for $\tilde I$ under a monomial order that eliminates $y_1,y_2$. For
each $P\in G$ set $P'=t^{-w(P)}P$ if $w(P)<0$ and $P'=\del_t^{w(P)}P$ if
$w(P)>0$ and let $G'=\{ P': P\in G\}$. Then
$G_0=G'\cap A_n[-\del_tt]$ generates $J\cap A_n[-\del_tt]$.
\pf
Note first that $G$ consists of $w$-homogeneous operators and so $w(P)$ is
well defined for $P\in G$.

Suppose $P\in G_0$. Hence $P\in\tilde I$. So $P=Q_{-1}\cdot (1-y_1y_2)+\sum
a_i\cdot (Q_i)^h$ where the $a_i$ are in $A_{n+1}[y_1,y_2]$. Since $P\in
A_n[-\del_tt]$, the 
substitution $y_i\to 1$ shows that $P=\sum a_i(1,1)\cdot
(Q_i)^h(1,1)=\sum a_i(1,1)\cdot Q_i\in J$. Therefore $G_0\sbe J\cap
A_n[-\del_tt]$. 

Now assume that $P\in J\cap A_n[-\del_tt]$. So $P$ is
$w$-homogeneous of weight 0. Also, $P\in J$ and $J$ is contained in
$I(1)$, the ideal of operators $Q(1)\sbe A_{n+1}$ for which $Q(y_1)\in
I$. By 
lemma \ref{oaku-lemma} below (taken from \cite{Oa}), $y_1^a
P\in I$ for some $a\in \Bbb  
N$. Therefore $P=(1-(y_1y_2)^a)P+(y_1y_2)^aP\in \tilde I$.

Let $G=\{P_1,\ldots,P_b,P_{b+1},\ldots,P_c\}$ and assume that $P_i\in
A_{n+1} $ if and only if $i\le b$. 
Buchberger algorithm gives a
standard expression $P=\sum a_iP_i$ with all $\ini(a_iP_i)\le \ini
(P)$. That implies that $a_{b+i}$ is zero for positive $i$ and $a_i$
does not contain $y_1,y_2$ for any $i$. 

Since $P, P_i$ are $w$-homogeneous,
the same is true for all $a_i$, from Buchberger algorithm. In fact,
$w(P)=w(a_i)+w(P_i)$ for all $i$. As $w(P)=0$ (and $t, \del_t$ are
the only variables with nonzero weight that may appear in $a_i$) we
find $a'_i\in A_n$ with $a_i=a_i'\cdot t^{-w(P_i)}$ or $a_i=a_i'\cdot
\del_t^{w(P_i)}$, 
depending on 
whether $w(P_i)$ is negative or positive. 

It follows that $P=\sum_1^b a_iP_i=\sum_1^b a_i'P_i'\in A_n[-\del_tt]\cdot
G_0$. \qed
\end{prop}
\begin{lem}
\mylabel{oaku-lemma}
Let $I$ be a $w$-homogeneous ideal in $A_{n+1}[y_1]$ with respect to the
weights introduced before the proposition and $I(1)$ defined as in the
proof of the proposition. Assume $P\in A_{n+1}$ is a $w$-homogeneous
operator. Then
$P\in I(1)$ implies 
$y_1^aP\in I$ for some $a$.
\pf
Note first that $y_1\to 1$ will not lead to cancellation of terms in any
homogeneous operator as $w(y_1)\not =0$.

If $P\in I(1)$, $P=\sum Q_i(1)$, with all $Q_i$ $w$-homogeneous in
$I$. Then the 
$y_1$-homogenization of $Q_i(1)$ will be a divisor of $Q_i$ and the
quotient will be some power of $y_1$, say
$y_1^{\eta_i}$. Homogenization of the equation 
$P=\sum Q_i(1)$ results in $y_1^{\eta}P=\sum Q_i(1)^h$ (since $P$ is
homogeneous) so that
\[\parbox{12.65cm}
{\hfill$y_{1}^{\eta+\max(\eta_i)}P=\sum y_1^{\max(\eta_i)-\eta_i}Q_i\in
I.$\hfill$ \Box$} \]
\end{lem}
So we have
\begin{alg}
\mylabel{ann-fs}
Input: $f\in R, L\sbe A_n$ such that $L$ is 
$f$-satuarated.

Output: Generators for $J^L(f^s)$.

Begin

\begin{enumerate}
\item For each generator $Q_i$ of $L$ compute the image $\phi(Q_i)$
under $x_i\to
x_i, t\to t-f, \del_i\to \del_i+f_i\del_t,\del_t\to\del_t$. Add $t-f$
to the list.

\item Homogenize all $\phi(Q_i)$ with respect to the new variable
$y_1$ relative to the weight $w$ introduced before proposition \ref{oaku}.

\item Compute a Gr\"obner basis for the ideal generated by 
$(\phi(Q_1))^h$, $\ldots$, $(\phi(Q_r))^h$, $1-y_1y_2$, $t-y_1f$ 
in $A_{n+1}[y_1,y_2]$
using an order that eliminates $y_1,y_2$.

\item Select the operators $\{ P_j\}_1^b$ in this basis which do not
contain $y_1, y_2$. 

\item For each $P_j$, $1\le j\le b$, if $w(P_j)>0$ replace $P_j$ by
$P_j'=\del_t^{w(P_j)}P_j$. Otherwise replace $P_j$ by
$P_j'=t^{-w(P_j)}P_j$. 

\item Return the new operators $\{P_j'\}_1^b$.
\end{enumerate}

End.
\end{alg}
In order to guarantee existence of the Bernstein polynomial $b^L_f(s)$
we assume for our next result that $L$ is holonomic.
\begin{cor}
\mylabel{b-poly}
Suppose $L$ is a holonomic ideal. 
If $J^L(f^s)$ is known or it is known that $L$ is $f$-saturated, then the
Bernstein polynomial $b_f^L(s)$ of $R_f\otimes_R 
A_n/L$ can be found from
 $(b^L_f(s))=A_n[s]\cdot(J^L(f^s),f)\cap K[s]$. 

Moreover, if $K\sbe \C$, 
suppose $b^L_f(s)=s^d+b_{d-1}s^{d-1}+\ldots+b_0$ and define 
$B=\max_{i}\{|b_i|^{1/(d-i)}\}$.
 In order to
find the smallest integer root of $b^L_f(s)$, one 
only needs to check the integers between $-2B$ and $2B$.

If in particular $L=(\del_1,\ldots,\del_n)$, it suffices to check the
integers between $-b_{d-1}$ and -1.
\pf 
If $L$ is $f$-saturated, propositions \ref{malgrange} and \ref{oaku}
enable us to find $J^L(f^s)$.  The first part follows then easily from the definition of $b^L_f(s)$:
as $(b_f^L(s)-P_f^L\cdot f)(f^s\otimes\bar 1)=0$ it is clear that
$b_f^L(s)$ is in $K[s]$ and in 
$A_n[s](J^L(f^s),f)$. Using an elimination order on $A_n[s]$,
$b^L_f(s)$ will be (up to a scalar 
factor) the unique element in the reduced Gr\"obner basis for
$J^L(f^s)+(f)$ that 
contains no $x_i$ nor $\del_i$.

Now suppose $K\sbe \C$, 
$|s|=2B\rho$ where $B$ is as defined above and $\rho>1$. 
Assume 
also that $s$ is a root of $b_f^L(s)$. We find
\begin{eqnarray}
(2B\rho)^d=|s|^d&=&|-\sum_0^{d-1}b_is^i|
%&\le&\sum_0^{d-1}|b_i|\cdot|s|^i \\
\le\sum_0^{d-1}B^{d-i}|s|^i\\
&=&B^d\sum_0^{d-1}(2\rho)^i
\le B^d((2\rho)^d-1),
\end{eqnarray}
using $\rho\geq 1$.
By contradiction, $s$ is not a root.

The
final claim is a consequence of Kashiwara's work \cite{K} where
it is proved that if $L=(\del_1,\ldots,\del_n)$ then  all roots of
$b_f^L(s)$ are negative and hence 
$-b_{n-1}$ is a 
lower bound for each single root.
\qed
\end{cor}

For purposes of reference we also list algorithms that compute the
Bernstein polynomial to a holonomic module and the localization of a
holonomic module.

\begin{alg}
\mylabel{b-poly-L}
Input: $f\in R, L\sbe A_n$ such that $A_n/L$ is holonomic and
$f$-torsionfree. 

Output: The Bernstein polynomial $b^L_f(s)$.

Begin
\begin{enumerate}
\item Determine $J^L(f^s)$ following algorithm \ref{ann-fs}. 

\item Find a reduced Gr\"obner basis for the ideal
$J^L(f^s)+A_n[s]\cdot f$ 
using an elimination order for $x$ and $\del$. 

\item Pick the unique element in that basis contained in $K[s]$ and
return it.
\end{enumerate}

End.
\end{alg}
\begin{alg}
\mylabel{D/L-loc-f}
Input: $f\in R, L\sbe A_n$ such that $A_n/L$ is holonomic and
$f$-torsionfree. 

Output: Generators for an ideal $J$ such that $R_f\otimes A_n/L\cong A_n/J$.

Begin
\begin{enumerate}
\item Determine $J^L(f^s)$ following algorithm \ref{ann-fs}. 
\item Find the Bernstein polynomial $b_f^L(s)$ using algorithm
\ref{b-poly-L}. 
\item Find the smallest integer root $a$ of $b_f^L(s)$ (using corollary
\ref{b-poly}, if $K\sbe \C$). 
\item Replace $s$ by $a$ in all generators for $J^L(f^s)$ and
return these generators.
\end{enumerate}

End.
\end{alg}
The algorithms \ref{ann-fs} and \ref{b-poly-L} appear already in
\cite{Oa} in the special case $L=(\del_1,\ldots,\del_n), A_n/L=R$. 

\section{Local cohomology as $A_n$-module}
\mylabel{sec-lc}
In this section we will combine the results from the previous sections to
obtain algorithms that compute for given $i,j,k\in \N, I\sbe  R$ the
local cohomology modules $H^k_I(R), H^i_\m(H^j_I(R))$ and the
invariants $\lambda_{i,j}(R/I)$ associated to $I$.
\subsection{Computation of $H^k_I(R)$}
\mylabel{subsec-lc}
Here we will describe an algorithm that takes in a finite
set of polynomials $\underline f=\{f_1,\ldots,f_r\}\subset R$ and
returns a 
presentation of $H^k_I(R)$ where $I=(f_1,\ldots,f_r)$. In particular,
if $H^k_I(R)$ is zero, then the algorithm will return the zero
presentation. 

Consider the  \v Cech complex associated to $f_1,\ldots,f_r$ in
$R$, 
\begin{equation}
\label{cechcomplex}
0\to R\to \bigoplus_1^r R_{f_i}\to \bigoplus_{1\le i<j\le r}R_{f_if_j}
\to\cdots\to R_{f_1\cdot\ldots\cdot f_r}\to 0.
\end{equation}
Its $k$-th cohomology group is the local cohomology module
$H^k_I(R)$.
The map 
\begin{equation}
\label{cechmap}
C^k=\bigoplus\limits_{1\le i_1<\cdots<i_k\le
r}R_{f_{i_1}\cdot\ldots\cdot f_{i_k}}\to \bigoplus\limits_{1\le
j_1<\cdots<j_{k+1}\le 
r}R_{f_{j_1}\cdot\ldots\cdot f_{j_{k+1}}}=C^{k+1}
\end{equation}
 is the sum of maps
\begin{equation}
\label{cechmap-parts}
R_{f_{i_1}\cdot\ldots\cdot f_{i_k}}\to R_{f_{j_1}\cdot\ldots\cdot
f_{j_{k+1}}}
\end{equation}
which are either zero (if $\{i_1,\ldots,i_k\}\not\subseteq
\{j_1,\ldots,j_{k+1}\}$) or send $\frac{1}{1}$ to
$\frac{1}{1}$, up to sign. 
Recall that $A_n/\Delta=
A_n/A_n\cdot(\del_1,\ldots,\del_n)\cong R$  and identify
$R_{f_{i_1}\cdot\ldots\cdot f_{i_k}}$ with
$A_n/J^\Delta((f_{i_1}\cdot\ldots\cdot f_{i_k})^s)|_{s=a}$ and
$R_{f_{j_1}\cdot\ldots\cdot f_{j_{k+1}}}$ with
$A_n/J^\Delta((f_{j_1}\cdot\ldots\cdot f_{j_{k+1}})^s)|_{s=b}$ where
$a,b$ are sufficiently small integers. By
proposition \ref{kashiwara} we may assume that $a=b\le 0$. Then the map
(\ref{cechmap}) is in the nonzero case multiplication from the right by
$(f_l)^{-a}$ where $l=\{j_1,\ldots,j_{k+1}\}\wo \{i_1,\ldots,i_k\}$,
again up to sign. 

It follows that the matrix representing the map $C^k\to C^{k+1}$ in
terms of $A_n$-modules is very easy to write down once the annihilator
ideals and Bernstein polynomials for all $k$- and $(k+1)$-fold products
of the $f_i$ are known: the entries are 0 or $\pm f_l^{-a}$ where
$f_l$ is the new factor. 

Let $\Theta^r_k$ be the set of $k$-element subsets of $1,\ldots,r$ and
for $\theta\in \Theta^r_k$ write $F_\theta$ for the product $\prod_{i\in
\Theta^r_k}f_{i}$.
We have demonstrated the correctness and finiteness of the following
algorithm.
\begin{alg}
\mylabel{lc-alg}

Input: $f_1,\ldots,f_r\in R; k\in \N$.

Output: $H_I^k(R)$ in terms of generators and relations as finitely
generated $A_n$-module.

Begin
\begin{enumerate}
\item Compute the annihilator ideal $J^\Delta((F_\theta)^s)$
and the Bernstein 
polynomial $b^\Delta_{F_\theta}(s)$ for all $(k-1)$-, $k$- and $(k+1)$-fold
products of 
${f_1}^s,\ldots,{f_r}^s$ as in \ref{ann-fs} and \ref{b-poly-L} (so
$\theta$ runs through $\Theta^r_{k-1}\cup \Theta^r_k\cup \Theta^r_{k+1}$).

\item Compute the smallest integer root $a_\theta$ for each
$b^\Delta_{F_\theta}(s)$, let $a$
be the minimum and replace $s$ by $a$ in all the annihilator ideals.

\item Compute the two matrices $M_{k-1},M_k$ representing the
$A_n$-linear maps 
$C^{k-1}\to C^k$ and $C^k\to C^{k+1}$ as explained in subsection
\ref{subsec-lc}. 

\item Compute a Gr\"obner basis $G$ for the kernel of the map
\[
\bigoplus_{\theta\in\Theta_k^r}A_n\to \bigoplus_{\theta\in\Theta^r_k} 
A_n/J^\Delta((F_\theta)^s)|_{s=a}\stackrel{M_k}{\longrightarrow}
\bigoplus_{\theta\in
\Theta^r_{k+1}}A_n/J^\Delta((F_\theta)^s)|_{s=a} \] as in \ref{kernel}.

\item Compute a Gr\"obner basis $G_0$ for the module 
\[
\im(M_{k-1})+\bigoplus_{\theta\in \Theta^r_k}
J^\Delta((F_\theta)^s)|_{s=a}\sbe \bigoplus_{\theta\in\Theta^r_k} 
A_n/J^\Delta((F_\theta)^s)|_{s=a}.
\]

\item Compute the remainders of all elements of $G$ with
respect to lifts of $G_0$ to $ \bigoplus_{\theta\in\Theta_k^r}A_n$. 

\item Return these remainders and $G_0$.
\end{enumerate}

End.
\end{alg}
The nonzero elements of $G$ generate the quotient $G/G_0\cong
H^k_I(R)$ so that
$H^k_I(R)=0$ if and only if all returned remainders are zero. 
\subsection{Computation of $H^i_\m( H^j_I(R))$}
As a second application of Gr\"ob\-ner basis computations over the
Weyl algebra  we
show now how to compute $H^i_\m (H^j_I(R))$.
Note that we cannot apply lemma \ref{malgrange} to $A_n/L=H^j_I(R)$
since $H^j_I(R)$ may well 
contain some torsion.

As in the previous sections, $C^j(R;f_1,\ldots,f_r)$ denotes the $j$-th
module in the 
\v Cech complex to $R$ and $\{f_1,\ldots,f_r\}$. 
Let $C^{\bullet\bullet}$ be the double complex with
$C^{i,j}=C^i(R;x_1,\ldots,x_n)\otimes_R C^j(R;f_1,\ldots,f_r)$, the
vertical maps $\phi^{\bullet\bullet}$ induced by the identity on the
first factor and the 
usual \v Cech maps on the second, whereas the horizontal maps
$\xi^{\bullet\bullet} $ are induced
by the \v Cech maps on the first factor and the identity on the
second. Since $C^i(R;x_1,\ldots,x_n)$ is $R$-projective, the column
co\-ho\-mo\-lo\-gy of $C^{\bullet\bullet}$ at $(i,j)$ is
$C^i(R;x_1,\ldots,x_n)\otimes_RH^j_I(R)$ and the induced horizontal maps
in the $j$-th row are
simply the maps in the \v Cech complex $C^\bullet(H^j_I(R);x_1,\ldots,x_n)$. 
It follows that
the row cohomology of the column cohomology at $(i_0,j_0)$ is
$H^{i_0}_\m(H^{j_0}_I(R))$. 

Now note that $C^{i,j}$ is a direct sum of modules $R_g$ where 
$g=x_{\alpha_1}\cdot\ldots\cdot x_{\alpha_i}\cdot
f_{\beta_1}\cdot\ldots\cdot f_{\beta_j}$. So the whole double complex
can be rewritten in terms of $A_n$-modules and $A_n$-linear maps using
\ref{D/L-loc-f}: 
\[
\diagram
{\,C^{i-1,j+1}\,}{\rto^{\,\,\xi^{i-1,j+1}}}&
	C^{i,j+1}\rto^{\xi^{i,j+1}}&
		C^{i+1,j+1}\\
C^{i-1,j}\rto^{\xi^{i-1,j}}\uto_{\phi^{i-1,j}}& 
	C^{i,j}\rto^{\xi^{i,j}}\uto_{\phi^{i,j}}& 
		C^{i+1,j}\uto_{\phi^{i+1,j}}\\
C^{i-1,j-1}\rto^{\xi^{i-1,j-1}}\uto_{\phi^{i-1,j-1}}&
	C^{i,j-1}\rto^{\xi^{i,j-1}}\uto_{\phi^{i,j-1}}& 
		C^{i+1,j-1}\uto_{\phi^{i+1,j-1}}
\enddiagram
\]
Using the comments in subsection \ref{double-kernel}, we may now
compute the modules $H^i_\m (H^j_I(R))$. More 
concisely, we have the following 
\begin{alg}
\mylabel{lclc-alg}
Input: $f_1,\ldots,f_r\in R; i_0,j_0\in \Bbb N$.

Output: $H^{i_0}_\m (H^{j_0}_I(R))$ in terms of generators and relations as
finitely generated $A_n$-module.

Begin.
\begin{enumerate}
\item For $i=i_0-1, i_0, i_0+1$ and $j=j_0-1,j_0,j_0+1$ compute the
annihilators $J^\Delta((F_\theta\cdot X_{\theta'})^s)$ and Bernstein
polynomials $b^\Delta_{F_\theta\cdot X_{\theta'}}(s)$ of $F_\theta\cdot
X_{\theta'}$ 
where $\theta \in \Theta^r_j, \theta'\in \Theta^n_i$ and $X_{\theta'}$
denotes in analogy to $F_\theta$ the product $\prod_{\alpha\in
\theta'}x_\alpha$.

\item Let $a$ be the minimum integer root of the product of all these
Bernstein polynomials and replace $s$ by $a$ in all the annihilators
computed in the previous step.

\item Compute the matrices to the $A_n$-linear maps
 $\phi^{i,j}:C^{i,j}\to 
C^{i,j+1}$ 
and $\xi^{i,j}:C^{i,j}\to C^{i+1,j}$, again for $i=i_0-1,i_0,i_0+1$ and
$j=j_0-1,j_0,j_0+1$. 

\item Compute Gr\"obner bases for the modules 
\[
G=\ker(\phi^{i_0,j_0})\cap
\left[ 
(\xi^{i_0,j_0})^{-1}(\im(\phi^{i_0+1,j_0-1}))\right]+\im(\phi^{i_0,j_0-1})
\] 
and 
$G_0=\xi^{i_0-1,j_0}(\ker(\phi^{i_0-1,j_0}))+\im(\phi^{i_0,j_0-1})$.

\item Compute the remainders of all elements of $G$ with
respect to $G_0$ 
and return these remainders together with $G_0$.
\end{enumerate}

End.
\end{alg}
The elements of $G$ will be generators for $H^{i_0}_\m
(H^{j_0}_I(R))$ and 
the elements of $G_0$ generate the relations that are not
syzygies.
\subsection{Computation of $\lambda_{i,n-j}(R/I)$}
In \cite{L-Dmod} it has been shown that $H^i_\m (H^j_I(R))$ is an injective
$\m$-torsion $R$-module of finite socle dimension $\lambda_{i,n-j}$
(which depends only on $i,j$ and $R/I$) and so
isomorphic to $(E_R(K))^{\lambda_{i,n-j}}$ where $E_R(K)$ is the
injective hull of $K$ over $R$. We
are now in a position that allows  computation of these invariants of $R/I$.
For, let $H^i_\m (H^j_I(R))$ be generated by $g_1,\ldots,g_l\in {A_n}^d$
modulo the relations $h_1,\ldots,h_e\in {A_n}^d$. Let $H$ be the module
generated by the $h_i$. We know that
$(A_n\cdot g_1+H)/H$ is $\m$-torsion and so it is possible (with trial
and error) to find a multiple of $g_1$, say $mg_1$ with $m$ a monomial
in $R$, such that $(A_n\cdot mg_1+H)/H$ is nonzero but $x_img_1\in H$ for
all $1\le i\le n$. Then the element $mg_1+H/H$ has annihilator equal
to $\m$ and hence generates an $A_n$-module isomorphic to
$A_n/A_n\cdot \m\cong E_R(K)$. The injection $A_n\cdot mg_1+H/H\into
A_n\cdot(g_1,\ldots,g_l)+H/H$ splits as map of $R$-modules because of 
injectivity and so the cokernel 
$A_n(g_1,\ldots,g_l)/A_n(mg_1,h_1,\ldots,h_e)$ is isomorphic to
$(E_R(K))^{\lambda_{i,n-j}-1}$. 

Reduction of the $g_i$ with respect to a Gr\"obner basis of the new
relation module and repetition of the previous will lead to
the determination of $\lambda_{i,n-j}$.
\subsection{Local cohomology in ambient spaces different from $\A^n_K$}
\mylabel{singular_spaces}
If $A$ equals $K[x_1,\ldots,x_n]$, $I\sbe A$, $X=\spec (A)$ and $V=\spec(A/I)$,
knowledge of $H^i_I(A)$ for all $i\in \N$ answers of course the
question about the local cohomological dimension of $V$ in $X$. It is
worth mentioning, that if $W\sbe X$ is a smooth variety containing $V$
then our algorithm \ref{lc-alg} for the computation of $H^i_I(A)$ also
leads to a determination of the local cohomological dimension of $V$
in $W$. Namely, if $J$ stands for the
defining ideal of $W$ in $X$ so that $R=A/J$ is the affine
coordinate ring of $W$  and if we set $c=\height(J)$, then it can be
shown that 
$H^{i-c}_{I}(R)=\hom_A(R,H^i_I(A))$ for all $i\in\N$.
As $H^i_I(A)$ is 
$I$-torsion (and hence $J$-torsion), $\hom_A(R,H^i_I(A))$ is zero if
and only if 
$H^i_I(A)=0$. It follows that the local cohomological dimension of $V$
in $W$ equals $\cd(A,I)-c$ and $\{q\in \N:H^q_I(A)\not =0\}=\{q\in
\N:H^{q-c}_I(R)\not =0\}$. 

If however $W$ is not smooth, no algorithms for the computation of
either $H^i_I(R)$ or $\cd(R,I)$ are known, irrespective of the
characteristic of the base field.

\section{Implementation and examples}
Some of the algorithms described above have been implemented as C-scripts
and tested on some examples. 
\subsection{}
The algorithm \ref{ann-fs} with $L=\Delta$ has been implemented by
Oaku  using the package Kan (see \cite{T}) which
is a postscript language for computations in the Weyl algebra and in
polynomial rings. An implementation for general $L$ is written by the current
author and part of a program that deals exclusively with computations
around local cohomology (\cite{W}). \cite{W} is theoretically able to compute
$H^i_I(R)$ for arbitrary $i, R={\Bbb Q}[x_1,\ldots,x_n], I\sbe R$ in the
above described terms of generators and relations, using algorithm
\ref{lc-alg}. It is expected that in the near future \cite{W} will
work for $R=K[x_1,\ldots,x_n]$ where $K$ is an arbitrary field of
characteristic zero and also algorithms
for computation of $H^i_\m (H^j_I(R))$ and $\lambda_{i,j}(R)$ will be
implemented, but see the comments in  
\ref{efficiency} below.
\begin{ex}
\label{example}
Let $I$ be the ideal in $R=K[x_1,\ldots,x_6]$ that is generated by the
$2\times 2$ minors $f,g,h$ of the matrix
$\left(\begin{array}{ccc}x_1&x_2&x_3\\x_4&x_5&x_6\end{array}\right)$.
Then $H_I^i(R)=0$ for $i<2$ and 
$H^2_I(R)\ne 0$ because $I$ is a height 2 prime and $H^i_I(R)=0$ for $i>3$
because $I$ is  
3-generated, so the only remaining case is $H^3_I(R)$. This module
 in
fact does 
not vanish, but until the discovery of our algorithm, its non-vanishing was a 
rather non-trivial fact. Our algorithm provides the first known proof of this 
fact by direct calculation.

Previously, Hochster pointed out that $H^3_I(R)$ is nonzero,
using the fact that the map $K[f,g,h]\to R$ splits (compare \cite{Hu-L},
Remark 3.13) and Bruns and Schw\"anzl (\cite{Br-S}, the corollary to
Lemma 2) provided a
topological proof of the nonvanishing of $H^3_I(R)$ via \'etale
cohomology.
 Both of these proofs are quite 
ingenious and work only in very special situations.

Using the
program \cite{W}, one finds that $H^3_I(R)$ is isomorphic to a cyclic
$A_6$-module 
generated by $1\in A_6$  subject to relations $x_1=\ldots =x_6=0$. 
This is a straightforward 
computational proof of the non-vanishing of $H^3_I(R)$. Of course this proof 
gives
more than simply the non-vanishing. Since the quotient of $A_6$ by the  
left ideal generated by $x_1,\dots,x_6$ is known to be isomorphic as an
$R$-module to  
$E_R(R/(x_1,\dots,x_6))$, the injective hull of $R/(x_1,\dots,x_6)=K$ in
the category  
of $R$-modules, our proof implies that 
$H^3_I(R)\cong E_R(K)$.
\end{ex}
\subsection{}
\mylabel{efficiency}
Computation of Gr\"obner bases in many variables is in general a time-
and space consuming enterprise. Already in (commutative) polynomial
rings the worst case performance for the number of elements in reduced
Gr\"obner bases 
is doubly exponential in the number of variables and the degrees of
the generators. In the (relatively small) 
example above $R$ is of dimension 6,
so that the intermediate ring $A_{n+1}[y_1,y_2]$ contains 16
variables. In view of these facts the following idea 
has proved useful. 

The general context in which lemma \ref{malgrange} and proposition
\ref{kashiwara} were stated allows successive localization of $R_{fg}$
in the following way. First one computes $R_f$ according to
algorithm \ref{D/L-loc-f} as quotient of $A_n$ by a certain holonomic 
ideal $L=J^\Delta(f^s)|_{s=a}, a\ll 0$.  
Then $R_{fg}$ may be
computed using  \ref{D/L-loc-f} again since $R_{fg}\cong R_g\otimes
A_n/L$. (Note that all
localizations of $R$ are automatically $f$-torsion free for $f\in R$
as $R$ is a domain.) This process
may be iterated for products with any finite number of factors. 
Note
also that the exponents for the various factors might be different. 
This requires some care as the following situations illustrate. Assume
first that $-1$ is the smallest integer root of the Bernstein polynomials
of $f$ and $g$ (both in $R$) with respect to the holonomic module
$R$. Assume further that $R_{fg}\cong A_n\cdot f^{-2}g^{-1}\supsetneq
A_n\cdot (fg)^{-1}$. Then $R_f\to R_{fg}$ can be written as
$A_n/\ann(f^{-1})\to A_n/\ann(f^{-2}\cdot g^{-1})$ sending $P\in A_n$ to
$P\cdot f\cdot g$. 

Suppose on the other hand that we are interested in $H^2_I(R)$ where
$I=(f,g,h)$ and we know that $R_f=A_n\cdot f^{-2}\supsetneq A_n\cdot
f^{-1}, R_g=A_n\cdot g^{-2}$ and $R_{fg}=A_n\cdot f^{-1}g^{-2}$. (In
fact, the degree 2 part of the \v Cech complex of example
\ref{example} consists of such localizations.) It is tempting to write
the embedding $R_f\to R_{fg}$ with the use of a Bernstein operator (if
$P_f(s) f^{s+1}=b^\Delta_f(s)f^s$ then take $s=-2$) but as $f^{-1}$ is not a
generator for $R_f$, $b^\Delta_f(-2)$ will be zero. In other words, we must
write $R_{fg}$ as $A_n/\ann((fg)^{-2})$ and then send $P\in
\ann(f^{-2})$ to $P\cdot g^2$. 

The two examples indicate how to write the \v Cech complex in terms
of generators and relations over $A_n$ while making sure that the maps
$C^k\to C^{k+1}$ can be made explicit using the $f_i$: the exponents
used in $C^i$ have to be at least as big as those in $C^{i-1}$ (for the
same $f_i$). 
\begin{rem}
\label{remark}
We suspect that for all holonomic $fg$-torsionfree 
modules $M=A_n/L$ 
we have  (with $R_g\otimes M\cong A_n/L'$):
\[\min\{s\in\Z:b_f^L(s)=0\}\le \min\{s\in\Z:b_f^{L'}(s)=0\}.\] 
This would have two nice consequences. 

First of all, it would
guarantee, that successive localization at the factors of a product
does not lead to matrices in the \v Cech complex with entries of
higher degree than localization at the product at once.

Secondly, if \ref{remark} were known to be true, we could proceed as
follows for the computation of $C^i(R;f_1,\ldots,f_r)$. First compute
$J^\Delta((f_i)^s)$ for all $i$, find all minimal integer 
Bernstein roots $\beta_i$
of $f_i$ on $R$ and substitute them into the appropriate annihilator
ideals. If from now on we want to use algorithm \ref{D/L-loc-f} in
order to compute $R_{f_{i_1}\cdot\ldots\cdot f_{i_k}\cdot
f_{i_{k+1}}}$ from $R_{f_{i_1}\cdot\ldots\cdot f_{i_k}}$ then we can
skip steps 2 and 3 of \ref{D/L-loc-f} as the remark gives us a lower bound for
the minimal integer Bernstein root of $f_{i_{k+1}}$ on
$R_{f_{i_1}\cdot\ldots\cdot f_{i_k}}$. (From the comments before
\ref{remark} it is also clear that we cannot hope to use a larger
value.) 
\end{rem}
The advantage of localizing $R_{fg}$ as $(R_f)_g$ is twofold. For
one, it allows the exponents of the various factors to be distinct
which is useful for the subsequent cohomology computation: it helps
to keep the degrees of the maps small. (So for example $R_{x\cdot
(x^2+y^2)}$ can be written as $A_n\cdot 
x^{-1} (x^2+y^2)^{-2}$ instead
of $A_n\cdot (x^{-2}\cdot (x^2+y^2)^{-2})$. 
On the other hand,  
since the computation of Gr\"obner bases is doubly exponential it
seems to be advantageous to break a big problem (localization at a
product) into many ``easy'' problems (successive localization).

An extreme case of this behaviour is our example \ref{example}: if we
compute $R_{fgh}$ as $((R_f)_g)_h$, the calculation uses
approximately 2.5 kB and lasts 32 seconds on a Sun
workstation using \cite{W}. If one tries to localize $R$ at the
product of the three generators at once, \cite{W} crashes after about
30 hours having used up the entire available memory (1.2 GB).

\bibliography{bib}
\bibliographystyle{abbrv}

\end{document}